\documentclass[reqno, 10pt]{amsart}

\usepackage{amsfonts}
\usepackage[latin1]{inputenc} 

\newtheorem{theorem}{Theorem}[section]

\newcommand{\N}{\mathbb{N}}
\newcommand{\Z}{\mathbb{Z}}

\makeatletter
\@addtoreset{equation}{section}
\makeatother

\renewcommand\thefigure{\thesection.\@arabic\c@figure}
\renewcommand\thetable{\thesection.\@arabic\c@table}

\bibdata{prob}
\bibstyle{alpha} 

\title[Spatial rumor Model]{A spatial stochastic model\\
for rumor transmission}

\author{Cristian F. Coletti, Pablo M. Rodríguez and Rinaldo B. Schinazi}

\date{}

\address{
\newline
UFABC - Centro de Matem\'atica, Computa\c{c}\~ao e Cogni\c{c}\~ao
\newline
Avenida dos Estados, 5001- Bangu - Santo André - São Paulo, Brasil
\newline
e-mail:  cristian.coletti@ufabc.edu.br
\newline
\newline
USP, Instituto de Ciências Matemáticas e de Computação
\newline  
Av. Trabalhador são-carlense 400 - Centro, CEP 13560-970, São Carlos, SP, Brasil
\newline
e-mail:  pablor@icmc.usp.br
\newline
\newline
Department of Mathematics, University of Colorado at Colorado Springs
\newline
Colorado Springs CO 80933, USA
\newline
e-mail: rschinaz@uccs.edu
}

\thanks{C.F.C. was partially supported by FAPESP (grant number 	09/52379-8); P.M.R. was supported by FAPESP (grant number 10/06967-2)}

\begin{document}
  
\maketitle

\begin{abstract}
We consider an interacting particle system representing the spread of a rumor by agents on the $d$-dimensional integer lattice. Each agent may be in any of the three states belonging to the set $\{0,1,2\}$. Here $0$ stands for ignorants, $1$ for spreaders and $2$ for stiflers. A spreader tells the rumor to any of its (nearest) ignorant neighbors at rate $\lambda$. At rate $\alpha$ a spreader becomes a stifler due to the action of other (nearest neighbor) spreaders. Finally, spreaders and stiflers forget the rumor at rate one. We study sufficient conditions under which the rumor either becomes extinct or survives with positive probability.
\end{abstract}

\section{Introduction} 
\label{sec:intro}
It is a well known fact the strong influence that a rumor may have on a social environment, for instance in fields such as politic, economy or social sciences. Although many mathematical models appeared in the scientific literature as an attempt to describe the behavior of this phenomenon we are, still, far from a complete understanding of its complexity.

Two classical stochastic models to describe the spread of a rumor were introduced by Daley and Kendall \cite{DK} and Maki and Thompson \cite{MT} for closed finite and homogeneously mixing populations. In both models the population is subdivided into three classes of individuals: ignorants, spreader and stiflers, and the rumor is propagated accordingly rules that depend of these classes. One of the main subjects under study in these processes is the survival or not of the rumor by analyzing the remaining proportion of ignorants when all spreaders have disappeared. 

After the first rigorous results, namely limit theorems for this proportion (Sudbury \cite{Sudbury} and Watson \cite{Watson}), many papers introduce modifications in the dynamic of the basic models in order to make them more realistic. Recent papers suggest generalizations which allow various contact interactions, the possibility of forgetting the rumor (Kawachi et al. \cite{Kawachi2}), long memory spreaders (Lebensztayn et al. \cite{RPRS2}), or a new class of ``uninterested'' individuals (Lebensztayn et al. \cite{RPRS}). However, all these models assume that the population is homogeneously mixed. 

Another issue to be considered is the role of space in the dynamic of rumor spreading. Recent articles deal with this question by considering a population on a graph. Machado et al. \cite{FMV} studied the behavior of discrete-time rumor processes when the population lives on $\Z$ and individuals have long range interactions. Also, there is increasing interest in understanding the diffusion of information on complex networks. The main tools used to study such processes are simulations and mean field approximations (see Isham et al. \cite{valerie} and references therein for more details).

The aim of this paper is to study a rumor process for a population on $\Z^d$, $d\geq 1$. We find sufficient conditions under which the rumor either becomes extinct or survives with positive probability. 
\section{The model and results} 
\label{sec:model}

The spatial rumor model we consider in this paper is a continuous-time Markov Process with state space $X=\{0,1,2\}^{\mathbb{Z}^d}$, i.e. at time $t$ the state of the process is some function $\eta_t: \mathbb{Z}^d \longrightarrow \{0,1,2\}$. We assume that at each site $x \in \mathbb{Z}^d$ there is an individual. An individual at $x \in \mathbb{Z}^d$ is said to be ignorant if $\eta(x)=0,$ spreader if $\eta(x)=1$ and stifler if $\eta(x)=2.$ If the system is in configuration $\eta \in X,$ the state of site $x$ changes from $i$ to $j$ at rate $c_{ij}(x,\eta).$ More precisely, the evolution at site $x$ is given by the transition rates

$$
\begin{array}{llllll}
c_{01}(x,\eta) &=& \lambda n_1(x,\eta)&  \ \ \ \ \ \ c_{10}(x,\eta) &=& 1  \\
\ & \ & \ & \ & \ & \ \\
c_{12}(x,\eta) &=& \alpha  n_1(x,\eta)& \ \ \ \ \ \ c_{20}(x,\eta) &=& 1 	
\end{array}
$$

\noindent where $n_1(x,\eta)=\displaystyle \sum_{||x-y||=1} 1\{\eta(y)=1\}$ is the number of nearest neighbors of site $x$ in state $1$ for the configuration $\eta .$ 

The main difference between this kind of model and the well known SIR for epidemic is due to the transition rate $c_{12}$. It represents the loss of interest of a spreader in transmitting the rumour after meeting another spreader.

\bigskip
\noindent {\it{Harris' graphical construction}}

Next we use Harris' graphical construction \cite{Harris} to construct the spatial rumor model.  

Consider a collection of independent Poisson processes $\{N^{x,y}_1 , N^{x,y}_2 , D^x : x,y \in \Z^d , \left\|x-y \right\|=1\}$. The processes $N^{x,y}_1 , N^{x,y}_2 , D^x$ have intensities $\lambda, \alpha$ and $1$, respectively. At each arrival time of $N^{x,y}_1$ if sites $x$ and $y$ are in states $1$ and $0$ respectively then, the state of site $y$ is updated to $1$. In a similar way, at each arrival time of $N^{x,y}_2$, if sites  $x$ and $y$ are both in state $1$ then $y$ changes to state $2$. Finally, at each arrival time of $D^x$ if there is a $1$ or a $2$ at site $x$ we replace it by a $0$. In this way we obtain a version of the spatial stochastic rumor model with the rates given above. In order to construct the process inside a finite space-time box it is enough to consider the Poisson arrival times inside that box. For further details on the graphical construction see Durrett \cite{durrett}.\\

\bigskip
\noindent {\it{Behavior of the rumor}}

Consider the spatial stochastic rumor model on $\Z^d$, for $d\geq 1$, and let $\lambda_c (d)$ be the critical value of the basic  $d$-dimensional contact process.

\begin{theorem} \label{1}
Let $\lambda\leq  \lambda_c (d)$. Then the rumor becomes extinct for all $\alpha \geq 0$.
\end{theorem}

Extinction has two different meanings in Theorem \ref{1} depending on whether the initial configuration has finitely many or infinitely many spreaders and stiflers. If the initial configuration has finitely many  spreaders and stiflers the rumor is said to become extinct if there is almost surely a finite random time after which all sites in $\mathbb{Z}^d$ are in state 0. If the initial configuration has infinitely many spreaders or stiflers the rumor is said to become extinct if for any fixed site there is almost surely a finite random time after which the site will stay in state 0 forever.

\begin{theorem} \label{2}
Let $\lambda > \lambda_c (d)$ fixed. Then there exists constants $0<\alpha_1 , \alpha_2<\infty$ such that
\begin{enumerate}
\item[(a)] if $\alpha < \alpha_1$ the rumor may survive;
\item[(b)] if $\alpha > \alpha_2$ the rumor becomes extinct.
\end{enumerate}
\end{theorem}
\section{The mean-field model} 
\label{sec:mean}

Let $u_i, i=0,1,2$ be the fraction of agents in state $i$. Note that $u_0+u_1+u_2=1$. Assuming independence between the states of nearest neighbor agents we arrive to the following coupled system of differential equations
\begin{eqnarray}
\frac{du_1}{dt}&=&\lambda u_1 (1-u_1-u_2) -  \alpha u_1^2-u_1, \nonumber \\
\frac{du_2}{dt}&=&\alpha u_1^2 - u_2 \nonumber
\end{eqnarray}
where $\lambda$ and $\alpha$ are the spreading and stifling rates respectively. Clearly, $(u_1,u_2)=(0,0)$ is a steady state for the system above. We call this the free rumor equilibrium state. Analogously to the theory of epidemic models we say that the spreading of the rumor is possible if this equilibrium is unstable. If $(0,0)$ is stable we say that the spreading of the rumor is not possible. The Jacobian matrix of the system at the steady state is
\[ J=\left( \begin{array}{cc}
\lambda-1 & 0 \\
0 & -1  \end{array} \right).\] 

The two eigenvalues are $\lambda-1$ and $-1$. Hence, the steady state $(0,0)$ is unstable if and only if $\lambda>1$. Unlike what happens for the spatial model, for the mean-field model whether the rumor spreads or not does not depend on $\alpha$.

\section{Proofs} 
\label{sec:proofs}

\noindent {\bf{Proof of Theorem \ref{1}}}

We consider a coupling between the spatial rumor process and the basic $d$-dimensional contact process $\xi_t$ with rates given by

\begin{eqnarray}\label{conpro}
c_{01}(x,\xi)= \lambda n_1(x,\xi)  \ \ \ \ \ \ c_{10}(x,\xi) = 1.  
\end{eqnarray}

At time $0$ we set $\xi_0(x)=1$ if $\eta_0(x)=1$ and $\xi_0(x)=0$ if $\eta_0(x)\not =1$. In words, we use as initial configuration for the contact process the initial configuration of the rumor process replacing all the $2$'s by $0$'s. We construct a version of the contact process  with the rates given above by using the same Poisson processes $N^{x,y}_1$ and $D^x$ considered in the graphical construction of the rumor process and ignoring the marks of the Poisson process $N^{x,y}_2$. 
More precisely, we use the rate $\lambda$ process $N^{x,y}_1$ to make appear a $1$ in $\xi$ and we use the rate $1$ process $D^x$ to make die a $1$ in $\xi$.  From this coupling it is not difficult to see that at all times the contact process has more $1$'s than the rumor process in the following sense. For any site $x$ in $\Z^d$ and any  time $t\geq 0$, if $\eta_t(x)=1$ then $\xi_t(x)=1$. Since we are assuming that $\lambda\leq \lambda_c(d)$ the contact process dies out for any initial configuration (finite or infinite). This implies that the spreaders in the rumor process die out. This in turn implies extinction of the rumor process.
The proof is complete. \\

\noindent {\bf{Proof of Theorem \ref{2}}}

\noindent {\it{Proof of (a)}}
We will compare the spatial rumor model with an oriented percolation model. In order for that we need some definitions. Consider $$\mathcal{L}_0=\{(m,n)\in \Z^2 : m+n \text{ is even}\}$$
and
$$
\begin{array}{rcl}
B = (-4L,4L)^d\times [0,T]&  \ \ \ \ \ \ & B_{m,n} = (2mL,nT)+B  \\
\ & \ &  \ \\
I = [-L,L]^d& \ \ \ \ \ \ &I_{m} = 2mL+I 	
\end{array}
$$
where $k=\sqrt L$ and $L$ and $T$ are values to be defined later.  We say that a site $(m,n)\in \mathcal{L}_0$ is open if and only if at time $nT$ there are no stiflers in $I_{m}$ and there are at least $k$ spreaders in $I_{m}$ and if at time $(n+1)T$ there are no stiflers in $I_{m-1}$ and $I_{m+1}$ and there are at least $k$ spreaders in each interval. Sites which are not open are called closed.
We will show that for any $\epsilon>0$, there exists a constant $\alpha_1 > 0$ such that for $\alpha < \alpha_1$ we get
\begin{equation}\label{open1}
\mathbb{P}[(m,n) \ \mbox{is open} \ ] \geq 1-5 \epsilon .
\end{equation}
By translation invariance, it is enough to show this for the site $(0,0)$. To simplify the notation we will suppose that $d=1$ in the proof. It is easy to check that all the arguments hold for any $d\geq 1$. We start by assuming $\alpha=0$ in $B$. This means that any spreader inside $B$ will not become stifler in $B$ except possibly at sites $-4L+1$ and $4L-1$. Stiflers are possible at  $-4L+1$ and $4L-1$ since we may have spreaders at sites $-4L$ and $4L$ (note that we assume nothing about the state of sites $-4L$ and $4L$ between times 0 and $T$). We assume that at time 0 there are no stiflers in $I$ but there could be stiflers elsewhere in $[-4L,4L]$ at time 0. However, the probability that all stiflers disappear by time $L$ in $[-4L+1,4L-1]$ is at least
$$1-(8L+1)\exp(-{L})\geq 1-\epsilon,$$
for $L$ large enough: there is no creation of stiflers and they die at rate 1.
Given that there are no stiflers in  $[-4L+1,4L-1]$ and that no stiflers can appear in 
$[-4L+1,4L-1]$ between times $L$ and $T$ the process of spreaders is a contact process with rates given by \eqref{conpro} in the finite space-time box $[-4L+1,4L-1]\times[L,T]$. 

Observe next that at time 0 there were at least $k$ spreaders and no stiflers in $I$ . Moreover, no stiflers can appear in $I$ between times $0$ and $T$. Hence, the spreaders in $I$ behave like a contact process and they survive at least as well as a contact process restricted to $I$. This is so because in our process spreaders could appear from outside $I$ into $I$ while this is not allowed in the contact process restricted to $I$. Now, the super-critical contact process restricted to the finite volume $\{1,2,\dots,k\}$ survives at least $e^k$ with probability $1-e^{-k}$, see Mountford \cite{Mountford}. It is not difficult to show that this in turn implies that the number of spreaders in $I$ at time $L$ is at least $M=\sqrt k$ with probability at least 
$$(1-e^{-M})^M\geq 1-\epsilon$$
for $L$ large enough.

So at time $L$ we have at least $M$ spreaders in $I$ and no stiflers in $[-4L+1,4L-1]$. 
Since $\lambda > \lambda_{c}(d)$, we can use well-known results of Bezuidenhout and Grimmett \cite{BG} for the supercritical contact process. In particular, given that a super-critical contact process does not die out the Shape Theorem (see Liggett \cite{Liggett}, p 128) ensures that the spreaders spread linearly. Hence, there is a constant $a>0$ such that by time $aL$ the spreaders have reached the sites $-4L$ and $4L$. Moreover, inside the cloud of spreaders there is a positive density of spreaders. Hence, by time $L+aL$ there are at least $k$ spreaders in $I_1$ and in $I_{-1}$ with probability at least $1-2\epsilon$ (one  $\epsilon$ takes care of the survival probability of $M$ spreaders and the other one of the Shape Theorem).  Set $T=L+aL$ where $L$ is large enough.

Now assume $\alpha>0$ and note that if we pick $\alpha$ small enough we can guarantee that the probability of having no marks of the process $N^{x,y}_2$ inside the finite space-time box $B$ is at least $1-\epsilon $. But this, implies that there exists a constant $\alpha_1>0$ such that if $\alpha <\alpha_1$ then
$$\mathbb{P}[(0,0) \ \mbox{is open} \ ] \geq 1-5\epsilon .$$

The rest of the proof is standard and depends on well known results on $K$-dependent oriented percolation (see Section 4 in Durrett \cite{durrett} for instance). Briefly, note that the event \{site $(m,n)\in \mathcal{L}_0$ is open\} only depends on the Poisson processes restricted to the box $B_{m,n}$. It is not difficult to see that there exist a constant $K$ such that if $||(m_i,n_i)-(m_j,n_j)||_{\infty}\geq K$, for $i\neq j$, the events  \{ site $(m_i,n_i)\in \mathcal{L}$ is open\} are independent. Since $\epsilon$ can be chosen small we conclude that there is percolation of open sites. Therefore, there will be spreaders at all times and this implies the survival of the rumor.\\   

\noindent {\it{Proof of (b)}}
The main idea for proving that the rumor becomes extinct is to compare the particle system $\eta_t$ to an oriented percolation process suitable defined on $\mathcal{L}=\mathbb{Z}^d \times \mathbb{Z}_+$.

Consider the following nested space-time regions:
\begin{eqnarray}
\Lambda_1 = [-2L,2L]^d \times [0,2T] \ \ \ \ \ \ \Lambda_2 = [-L,L]^d \times [T,2T].
\end{eqnarray}
Let $\Delta$ be the boundary of $\Lambda_1$:
\[
\Delta = \{(x,t) \in \Lambda_1 : |x_i|=2L \ \mbox{for some} \ i, i=1, \ldots d \ \mbox{or} \ t=0\}.
\]
We associate to each lattice point $(x,t)$ a $0-1$ valued random variable $\omega(x,t)$. If $\omega(x,t)=1$ we say that site $(x,t)$ is open, otherwise we say that site $(x,t)$ is closed. We say that a site $(x,t) \in \mathcal{L}$ is open if and only if for the process restricted to $\Lambda_1 + (x,t)$ the box $\Lambda_2 + (x,t)$ is empty regardless of the states of sites in the boundary $\Delta + (x,t)$. Sites which are not open are called closed.

We will show that for any $\epsilon > 0$, there exists $\alpha_2 >0$ such that if $\alpha > \alpha_2$ then:
\begin{equation}\label{open}
\mathbb{P}[(x,t) \ \mbox{is open} \ ] \geq 1-2\epsilon .
\end{equation}

Note that by translation-invariance it suffices to consider the site $(0,0) \in \mathcal{L}$. Suppose that $(0,0)$ is closed. That is, there is a spreader or stifler inside $\Lambda_2$. In order to reach $\Lambda_2$ the rumor must have originated from a site on the boundary $\Delta$. Moreover, there must be a path of nearest neighbor sites inside  $\Lambda_1$ that spread the rumor from  some site in $\Delta$ to some site in  $\Lambda_2$. Now,  since $\Lambda_1$ is a finite space-time box we can pick $\alpha$ so large that the following event $A$ happens with probability at least $1-\epsilon$. Let $A$ be the event that every time a spreader is born in $\Lambda_1$, before anything else happens, the new spreader becomes a stifler or the nearest neighbor which gave birth to the new spreader becomes a stifler. Hence, given the event $A$ any path of sites that spread the rumor from $\Delta$ to $\Lambda_2$ has only one spreader at any given time. Recall that spreaders die at rate 1. Moreover,  a path of spreaders either originates at the bottom of $\Delta$ or on the sides of $\Delta$. In the first case the path of spreaders must survive for at least $T$. In the second case the path of spreaders must  have covered a distance at least $L$ and therefore must have survived at least $bL$ for some constant $b>0$.
Hence,
\begin{equation}
\mathbb{P}[(0,0) \ \mbox{is open} |A ] \geq 1-\left(4L+1\right)^d e^{-T}-(2d)(2T)(4L+1)^{d-1}e^{-bL}.
\end{equation}

By taking $L=T$ large enough we get
\begin{equation}
\mathbb{P}[(0,0) \ \mbox{is open} |A ] \geq 1-\epsilon
\end{equation}
 For $\alpha$ large enough we have
 $$P(A)>1-\epsilon.$$
 Therefore, we get that for  $\alpha$ large enough
 $$\mathbb{P}[(0,0) \ \mbox{is open} ] \geq 1-2\epsilon.$$

By translation invariance, the same is true for any site $(x,t)$ in $\mathcal{L}$.

We now make $\mathcal{L}$ into a graph. Let $\Lambda (x,t) = \Lambda_1 + (x,t)$. Then, draw an oriented edge from $(x,t)$ to $(y,t^{\prime})$ if and only if $t \leq t^{\prime}$ and $\Lambda (x,t) \cap \Lambda (y,t^{\prime}) \neq \emptyset$. The open sites in the resulting directed graph defines a percolation model as follows. We say that $(y,t^{\prime})$ can be reached from $(x,t)$ and write $(x,t) \rightarrow (y,t^{\prime})$ if there is a sequence of sites $x_0=x, \ldots, x_n=y$ and instant times $t_0=t, \ldots, t_n=t^{\prime}$ such that: first, there is an oriented edge from $(x_k,t_k)$ to $(x_{k+1},t_{k+1})$ for $0 \leq k < n$; second, $\omega (x_k,t_k)=1$ for $0 \leq k \leq n.$ Note that the state of sites are not independent. However, they are $K$ dependent in the following sense; there exists a constant $K$ depending only on the dimension $d$ such that if the distance between sites $(x,t)$ and $(y,t^{\prime})$ is larger than $K$ then the state of the corresponding sites are independent.

The rest of the proof is somewhat standard. We refer the reader to Van Den Berg et al. \cite{bergs} for more details. The crucial point is that a site with a spreader  or a stifler in the rumor process corresponds with a path of closed sites in the percolation process. By taking $\epsilon$ small enough it is possible to make the probability of a path of closed sites decay exponentially fast with its length. This in turn implies that for any fixed site in the rumor process there is a finite random time after which the site is in state 0. This completes the proof.


\end{document}